\documentclass{amsart}
\usepackage{graphicx}
\usepackage{amssymb}
\usepackage{epstopdf}
\usepackage{nicefrac}
\usepackage{enumerate}
\usepackage{hyperref}
\usepackage{float}
\usepackage{color}
\usepackage{pst-all}
\usepackage{tabularx}
\usepackage{soul}

\newtheorem{thm}{THEOREM}[section]

\newtheorem{cor}[thm]{COROLLARY}

\newtheorem{defn}[thm]{DEFINITION}
\newtheorem{ex}[thm]{EXAMPLE}

\newtheorem{lemma}[thm]{LEMMA}

\newtheorem{prop}[thm]{PROPOSITION}

\newtheorem{remark}[thm]{REMARK}

\newcommand{\ds}{\displaystyle}

\newcommand{\cH}{{\mathcal H}}

\newcommand{\cO}{{\mathcal O}}

\newcommand{\CO}{{\rm CO}} 

\newcommand{\cU}{{\mathcal U}}

\newcommand{\cW}{{\mathcal W}}
\newcommand{\diam}{{\rm diam}} 
\newcommand{\dist}{{\rm dist}} 
\newcommand{\dX}{d_{\fX}} 
\newcommand{\e}{{\varepsilon}} 

\newcommand{\Fix}{{\rm Fix}} 

\newcommand{\fX}{{\mathfrak{X}}}

\newcommand{\G}{\Gamma}

\newcommand{\Homeo}{{\rm Homeo}} 
 
\newcommand{\Id}{{\rm Id}} 
\newcommand{\Iso}{{\rm Iso}} 
\newcommand{\mZ}{{\mathbb Z}}

\newcommand{\Perm}{{\rm Perm}} 
\newcommand{\vc}{{\rm vc}}

\newcommand{\whg}{\widehat{g}}

\newcommand{\whGamma}{\widehat{\Gamma}}

\newcommand{\whrho}{{\widehat{\rho}}}

\newcommand{\cR}{{\mathcal R}}

 \newcommand{\whalpha}{\widehat{\alpha}}
\newcommand{\whtheta}{{\widehat{\theta}}}
\newcommand{\cZ}{{\mathcal Z}}

\begin{document}

\title{Nilpotent   Cantor actions}

\author{Steven Hurder}
\address{Steven Hurder, Department of Mathematics, University of Illinois at Chicago, 322 SEO (m/c 249), 851 S. Morgan Street, Chicago, IL 60607-7045}
\email{hurder@uic.edu}

\author{Olga Lukina}
\address{Olga Lukina, Faculty of Mathematics, University of Vienna, Oskar-Morgenstern-Platz 1, 1090 Vienna, Austria}
\email{olga.lukina@univie.ac.at}

\thanks{Version date: November 12, 2020; rev. March 22, 2021}

\thanks{2010 {\it Mathematics Subject Classification}. Primary:  37B05, 37C15, 37C85; Secondary: 57S10}

\thanks{OL is supported by FWF Project P31950-N35}

  \thanks{Keywords: minimal   Cantor actions,  topological orbit equivalence,  return equivalence,   topologically free actions}

  \begin{abstract}
  A  nilpotent   Cantor action  is a  minimal equicontinuous   action $\Phi \colon \G \times \fX \to \fX$ on a Cantor space $\fX$,  where $\G$ contains a finitely-generated nilpotent subgroup $\G_0 \subset \G$ of finite index. 
  In this note, we show that these actions are distinguished among general Cantor actions: any effective action of a finitely generated group on a Cantor space, which is continuously orbit equivalent to a nilpotent Cantor action, must itself be a nilpotent Cantor action. As an application of this result, we obtain     new invariants of nilpotent Cantor actions under continuous orbit equivalence.
     \end{abstract}

\maketitle

 \vspace{-.1in}


\section{Introduction}\label{sec-intro}
  
  Let $\G$ be a countably generated   group, and let $\Phi \colon \G \times \fX \to \fX$, also denoted by $(\fX,\G,\Phi)$,  be   an action of $\G$ on a topological space $\fX$.   We say it   is  a \emph{Cantor action} if $\fX$ is a Cantor space.

A \emph{nilpotent   Cantor action} is a  minimal equicontinuous Cantor action  $(\fX,\G,\Phi)$,  where $\G$ contains a finitely-generated nilpotent subgroup $\G_0 \subset \G$ of finite index. Nilpotent Cantor actions arise in a variety of contexts, which motivates our interest in this  class of actions. 

A  minimal equicontinuous   Cantor action  is   called a \emph{generalized odometer} in the works \cite{CortezPetite2008,CortezMedynets2016,GPS2019,Li2018}, and when $\G = \mZ$  then $(\fX,\G,\Phi)$ is just a classical  odometer, which has been extensively studied   \cite{Downarowicz2005}.  In this work we study properties of generalized odometers given by a virtually nilpotent group action.

A classical odometer  is determined up to topological conjugacy by a supernatural number associated to the action (see Bing \cite{Bing1960}, Aarts and Fokkink \cite{AartsFokkink1991}). When $\G$ is a finitely-generated free abelian group,  then the generalized odometers are   classified up to continuous orbit equivalence in the works by Cortez and Medynets \cite{CortezMedynets2016} and Giordano, Putnam and Skau \cite{GPS2019}.
The nilpotent Cantor actions can   be considered as  the ``simplest'' class of Cantor actions whose classification problem is unresolved. One goal of their study is   to find invariants of the actions which distinguish them  up to topological conjugacy, or better, up to  continuous orbit equivalence.
  
 Another motivation for studying nilpotent Cantor actions is that they arise in the classification of \emph{renormalizable groups}; that is, finitely generated groups  which admit a proper self-embedding whose image has finite index   \cite{HLvL2020}. The works by Cornulier \cite{Cornulier2016} and Der\'e    \cite{Dere2017} give general criteria for when a nilpotent group admits such a self-embedding.  Renormalizable groups arise in a number of geometric and dynamical contexts, such as in the study of laminations with the shape of a compact manifold \cite{CHL2020}, and in the classification of generalized Hirsch foliations   \cite{BHS2006}.

 There is a well-developed theory of the ergodic properties of measure-preserving ergodic actions of nilpotent groups (for example, see the book by Host and Kra \cite{HostKra2018}), but not so much for the topological dynamics of nilpotent Cantor actions. This paper makes a contribution to their study. The terms in the following result are defined  in   Section~\ref{sec-basics}.

   \begin{thm}\label{thm-nilconjugate}
   Let     $(\fX_1,\G_1,\Phi_1)$ be a nilpotent Cantor action which is   continuously orbit equivalent to a   Cantor action $(\fX_2,\G_2,\Phi_2)$, then the actions $\Phi_1$ and $\Phi_2$ are return equivalent. Moreover, if both actions are effective, or faithful, then $(\fX_2,\G_2,\Phi_2)$ is  a nilpotent Cantor action.  If both actions   are topologically free, then $\G_1$ and $\G_2$ have nilpotent subgroups of finite index which are isomorphic, and so in particular,  $\G_1$ and $\G_2$  are commensurable.
       \end{thm}

Given Cantor actions  $(\fX_1,\G_1,\Phi_1)$ and $(\fX_2,\G_2,\Phi_2)$, we can replace them with effective actions by considering the actions of the quotient groups $\G_1' = \G_1/\ker(\Phi_1)$ and $\G_2' = \G_2/\ker(\Phi_2)$ to which Theorem~\ref{thm-nilconjugate} applies.

Example~\ref{ex-comm} shows that the conclusion that $\G_2$ contains a nilpotent subgroup of finite index is best possible.
 Example~\ref{ex-notcomm} shows that if the actions are not topologically free, then the finite-index nilpotent subgroups of $\G_1$ and $\G_2$ need not be isomorphic, or even commensurable. 

  Theorem \ref{thm-nilconjugate}   suggests  that invariants of continuous orbit equivalence for nilpotent Cantor actions must be ``virtual'' in nature, and depend on properties of nilpotent groups of a special nature. Here is one such invariant.

  The \emph{virtual nilpotency class} $\textrm{vc}(\G)$ of a finitely-generated virtually nilpotent group $\G$ is defined as the length of a central series for a torsion-free nilpotent subgroup of finite index. This is discussed further in Section~\ref{sec-class}. 
 The proof of Theorem~\ref{thm-nilconjugate} yields the following result.

  \begin{thm}\label{thm-vc}
   Let     $(\fX_1,\G_1,\Phi_1)$ and  $(\fX_2,\G_2,\Phi_2)$ be effective    Cantor actions, with $\G_1$ and $\G_2$ finitely generated.
 Suppose that     $(\fX_1,\G_1,\Phi_1)$ is a nilpotent action, and    $(\fX_2,\G_2,\Phi_2)$ is       continuously orbit equivalent to $(\fX_1,\G_1,\Phi_1)$. Then $(\fX_2,\G_2,\Phi_2)$ is a nilpotent Cantor action, and $\textrm{vc}(\G_1) = \textrm{vc}(\G_2)$.  
        \end{thm}

 As a second  application, in the work \cite{HL2020b} the authors study the \emph{asymptotic prime spectrum} of an equicontinuous Cantor action, which is a generalization of the invariant which classifies equicontinuous actions of $\mZ$ as in \cite{AartsFokkink1991,Bing1960}. Theorem~\ref{thm-nilconjugate} implies that the asymptotic prime spectrum is an invariant of nilpotent Cantor actions under continuous orbit equivalence.

Associated to an equicontinuous Cantor action $(\fX,\G,\Phi)$ is a reduced $C^*$-algebra $C_r^*(\fX,\G,\Phi)$ with a natural choice of Cartan subalgebra, as defined by   Renault \cite{Renault2008}.
  Renault  studies the properties of Cartan subalgebras and their relation to dynamical systems. 
The results of \cite{Renault2008} and the structure theory for   $C^*$-algebras of Type I, as in Arveson \cite{Arveson1976}, can be used to define invariants of nilpotent Cantor actions, which are invariants under continuous orbit equivalence by   Theorem~\ref{thm-nilconjugate}.

In Section~\ref{sec-basics}   we explain the terminology   and recall necessary preliminary results for the  proof of Theorem~\ref{thm-nilconjugate}. In Section~\ref{sec-equicontinuous} we show that equicontinuity is preserved by continuous orbit equivalence, and in Section~\ref{sec-return} we give a result showing that then the Cantor actions  become return equivalent.
The proof of Theorem~\ref{thm-nilconjugate} is then given in Section~\ref{sec-nilpotent}.    The virtual nilpotent class    of a virtually nilpotent group and nilpotent Cantor action are defined in Section~\ref{sec-class}, where we   give a proof of Theorem~\ref{thm-vc}.

 \section{Cantor actions}\label{sec-basics}

We recall some of the basic 
 properties of     Cantor actions.
References for the results described below are     the text by Auslander \cite{Auslander1988}, the papers by Cortez and Petite  \cite{CortezPetite2008}, Cortez and Medynets  \cite{CortezMedynets2016},    and   the authors' works  \cite{DHL2016c} and   \cite[Section~3]{HL2019b}.

\subsection{Basic concepts}\label{subsec-basics}

Let  $(\fX,\G,\Phi)$   denote an action  $\Phi \colon \G \times \fX \to \fX$. We   write $g\cdot x$ for $\Phi(g)(x)$ when appropriate.
The orbit of  $x \in \fX$ is the subset $\cO(x) = \{g \cdot x \mid g \in \G\}$. 
The action is \emph{minimal} if  for all $x \in \fX$, its   orbit $\cO(x)$ is dense in $\fX$.

   Let $N(\Phi) \subset \G$ denote the kernel of the action homomorphism $\Phi \colon \G \to \Homeo(\fX)$. The action is said to be \emph{effective} if $N(\Phi)$ is the trivial group. That is, the homomorphism $\Phi$ is faithful, and one also says that the action  is faithful.
    
 An action  $(\fX,\G,\Phi)$ is \emph{equicontinuous} with respect to a metric $\dX$ on $\fX$, if for all $\e >0$ there exists $\delta > 0$, such that for all $x , y \in \fX$ and $g \in \G$ we have  that 
 $\ds  \dX(x,y) < \delta$ implies   $\dX(g \cdot x, g \cdot y) < \e$.
The property of being equicontinuous    is independent of the choice of the metric   on $\fX$,  compatible with the topology of $\fX$.

\vfill
\eject

Now assume that $\fX$ is a Cantor space. 
Let $\CO(\fX)$ denote the collection  of all clopen (closed and open) subsets of  $\fX$, which forms a basis for the topology of $\fX$. 
For $\phi \in \Homeo(\fX)$ and    $U \in \CO(\fX)$, the image $\phi(U) \in \CO(\fX)$.  
The following   result is folklore, and a proof is given in \cite[Proposition~3.1]{HL2018b}.
 \begin{prop}\label{prop-CO}
 For $\fX$ a Cantor space, a minimal   action   $\Phi \colon \G \times \fX \to \fX$  is  equicontinuous  if and only if  the $\G$-orbit of every $U \in \CO(\fX)$ is finite for the induced action $\Phi_* \colon \G \times \CO(\fX) \to \CO(\fX)$.
\end{prop}

We say that $U \subset \fX$  is \emph{adapted} to the action $(\fX,\G,\Phi)$ if $U$ is a   \emph{non-empty clopen} subset, and for any $g \in \G$, 
if $\Phi(g)(U) \cap U \ne \emptyset$ implies that  $\Phi(g)(U) = U$.     Given  $x \in \fX$ and clopen set $x \in W$, there is an adapted clopen set $U$ with $x \in U \subset W$. (For a proof of this, see  \cite[Proposition~3.1]{HL2018b}.)
It follows that the adapted sets containing a point $x \in \fX$ form a local base for the topology. We single out a choice of a base with the following definition:

\begin{defn}\label{def-adaptednbhds}
Let  $(\fX,\G,\Phi)$   be a Cantor    action.
A properly descending chain of clopen sets $\cU = \{U_{\ell} \subset \fX  \mid \ell \geq 0\}$ is said to be an \emph{adapted neighborhood basis} at $x \in \fX$ for the action $\Phi$,   if
    $x \in U_{\ell +1} \subset U_{\ell}$ for all $ \ell \geq 0$ with     $\cap_{\ell > 0}  \ U_{\ell} = \{x\}$, and  each $U_{\ell}$ is adapted to the action $\Phi$.
\end{defn}

A key property    is that   for  $U$   adapted,   the set of ``return times'' to $U$, 
 \begin{equation}\label{eq-adapted}
\G_U = \left\{g \in \G \mid g \cdot U  \cap U \ne \emptyset  \right\}  
\end{equation}
is a subgroup of   $\G$, called the \emph{stabilizer} of $U$.      
  Then for $g, g' \in \G$ with $g \cdot U \cap g' \cdot U \ne \emptyset$ we have $g^{-1} \, g' \cdot U = U$, hence $g^{-1} \, g' \in \G_U$. Thus,  the  translates $\{ g \cdot U \mid g \in \G\}$ form a finite clopen partition of $\fX$, and are in 1-1 correspondence with the quotient space $X_U = \G/\G_U$. Then $\G$ acts by permutations of the finite set $X_U$ and so the stabilizer group $\G_U \subset\G$ has finite index.  Note that this implies that if $V \subset U$ is a proper inclusion of adapted sets, then the inclusion $\G_V \subset \G_U$ is also proper.

\subsection{Fixed points for Cantor actions}\label{subsec-regularity}
 We next consider  the structure of the sets of fixed points for a Cantor action $(\fX,\G,\Phi)$.

  The  action  is \emph{free} if for all $x \in \fX$ and $g \in \G$,   $g \cdot x = x$ implies that $g = e$,   the identity of the group. The \emph{isotropy group} of $x \in \fX$ is the subgroup 
\begin{equation}\label{eq-isotropyx}
\G_x = \{ g \in \G \mid g \cdot x = x\} \ . 
\end{equation}

Let $\Fix(g) = \{x \in \fX \mid g \cdot x = x \}$, and introduce the \emph{isotropy set}
\begin{equation}\label{eq-isotropy}
 \Iso(\Phi) = \{ x \in \fX \mid \exists ~ g \in \G ~ , ~ g \ne id ~, ~g \cdot x = x    \} = \bigcup_{e \ne g \in \G} \ \Fix(g) \ . 
\end{equation}

  \begin{defn}\cite{BoyleTomiyama1998,Li2018,Renault2008} \label{def-topfree}
  $(\fX,\G,\Phi)$ is said to be \emph{topologically free}  if  $\Iso(\Phi) $ is meager in $\fX$. 
 \end{defn}
 
Note that if $\Iso(\Phi)$ is meager, then $\Iso(\Phi)$ has empty interior. That is, if there exists a non-identity element $g \in \G$ such that $\Fix(g)$ has interior, then the action is not topologically free.
 
   The   notion of topologically free actions was introduced by Boyle in his thesis \cite{Boyle1983}, and later  used  in the works by Boyle and Tomiyama \cite{BoyleTomiyama1998}   for the study of classification of general Cantor actions,   by   Renault \cite{Renault2008}
     for the study of the $C^*$-algebras associated to Cantor actions, and by   Li \cite{Li2018} for proving    rigidity properties of Cantor actions.

The notion of a 
  \emph{quasi-analytic} action, which   was introduced in the  works  of    {\'A}lvarez L{\'o}pez,   Candel, and Moreira Galicia  \cite{ALC2009,ALM2016}, is an alternative formulation of the topologically free property which generalizes to group actions   where the acting group can be countable or profinite. 
  \begin{defn}\label{def-qa}
 An action $\Phi \colon H \times \fX \to \fX$, where 
    $H$ is a topological group and  $\fX$ a Cantor space,    is   \emph{quasi-analytic} if for each clopen set $U \subset \fX$, 
  if  the action of $g \in H$ satisfies $\Phi(g)(U) = U$ and the restriction $\Phi(g) | U$ is the identity map on $U$, 
  then $\Phi(g)$ acts as the identity on $\fX$.  
  \end{defn}

A topologically free action, as in Definition~\ref{def-topfree}, is quasi-analytic. That is, the isotropy set   \eqref{eq-isotropy}     has non-empty interior if $(\fX, H,  \Phi)$  is not quasi-analytic. Conversely,  the Baire Category Theorem implies that a quasi-analytic   effective action of a countable group  is topologically free \cite[Section~3]{Renault2008}.

  A local formulation of the quasi-analytic property was introduced  in the works \cite{DHL2016c,HL2018a}, and has proved very useful for the study of the dynamical properties of Cantor actions. 
    \begin{defn}  \label{def-LQA}  
   An action $\Phi \colon H \times \fX \to \fX$, where 
    $H$ is a topological group and  $\fX$ a Cantor metric space with metric $\dX$,   is   \emph{locally quasi-analytic}  if there exists $\e > 0$ such that for any non-empty open set $U \subset \fX$ with $\diam (U) < \e$,  and  for any non-empty open subset $V \subset U$,  if the action of $g \in H$ satisfies $\Phi(g)(V) = V$ and the restriction $\Phi(g) | V$ is the identity map on $V$,    then $\Phi(g)$ acts as the identity on   $U$.  
\end{defn}

\subsection{Equivalence of Cantor actions}\label{subsec-equivalence}

We recall three notions of equivalence of  Cantor actions which we use in this work. This first and strongest notion is the following, as used   in   \cite{CortezMedynets2016,HL2018b,Li2018}: 

 \begin{defn} \label{def-isomorphism}
Cantor actions $(\fX_1, \G_1, \Phi_1)$ and $(\fX_2, \G_2, \Phi_2)$ are said to be \emph{isomorphic}  if there is a homeomorphism $h \colon \fX_1 \to \fX_2$ and group isomorphism $\Theta \colon \G_1 \to \G_2$ so that 
\begin{equation}\label{eq-isomorphism}
\Phi_1(g) = h^{-1} \circ \Phi_2(\Theta(g)) \circ h   \in   \Homeo(\fX_1) \   \textrm{for  all} \ g \in \G_1 \ .
\end{equation}
 \end{defn}

The notion of \emph{return equivalence} for Cantor actions is  defined next. This equivalence is  weaker than the notion of isomorphism, and is natural when considering the Cantor systems defined by the holonomy actions for matchbox manifolds, as considered in   the works  \cite{HL2018a,HL2018b,HL2019b}.

For a minimal equicontinuous Cantor action $(\fX, \G, \Phi)$ and   an adapted set $U \subset \fX$, by a small abuse of  notation, we use $\Phi_U$ to denote both the restricted action $\Phi_U \colon \G_U \times U \to U$ and the induced quotient action $\Phi_U \colon H_U \times U \to U$ for $H_U = \Phi(\G_U) \subset \Homeo(U)$. Then $(U, H_U, \Phi_U)$ is called the \emph{holonomy action} for $\Phi$, in analogy with the case where $U$ is a transversal to a matchbox manifold, and 
  $H_U$ is the holonomy group for this transversal.

   \begin{defn}\label{def-return}
Two minimal equicontinuous Cantor actions $(\fX_1, \G_1, \Phi_1)$ and $(\fX_2, \G_2, \Phi_2)$ are said to be  \emph{return equivalent} if there exists 
  an adapted set $U \subset \fX_1$ for the action $\Phi_1$   and  
  an adapted set $V \subset \fX_2$ for the action $\Phi_2$,
such that   the  restricted actions $(U, H_{1,U}, \Phi_{1,U})$ and $(V, H_{2,V}, \Phi_{2,V})$ are isomorphic.
\end{defn}

The  notion of \emph{continuous orbit equivalence} for   Cantor actions  was introduced in     \cite{Boyle1983,BoyleTomiyama1998}, and plays a fundamental role in various approaches to the classification of these actions \cite{Renault2008}.
Consider the equivalence relation on $\fX$ defined by an action    $(\fX,\G,\Phi)$,
 \begin{equation}\label{eq-ERX}
\cR(\fX, \G, \Phi) \equiv \{(x,  \gamma x)) \mid x \in \fX, \gamma \in \G\} \subset \fX \times \fX ~.
\end{equation}
Given   actions $(\fX_1,\G_1,\Phi_1)$ and $(\fX_2, \G_2, \Phi_2)$, we say they are \emph{orbit equivalent} if there exist a bijection   $h \colon \fX_1 \to \fX_2$ which maps 
$\cR(\fX_1, \G_1, \Phi_1)$   onto   $\cR(\fX_2, \G_2, \Phi_2)$, and similarly for the inverse map $h^{-1}$.

\begin{defn}\label{def-torb1}
Let $(\fX_1,\G_1,\Phi_1)$ and $(\fX_2, \G_2, \Phi_2)$ be      Cantor actions.
A \emph{continuous orbit equivalence}  between the actions  is a homeomorphism    $h \colon \fX_1 \to \fX_2$ which is an orbit equivalence, and there exist continuous functions $\alpha \colon \G_1 \times \fX_1 \to \G_2$ and $\beta \colon \G_2 \times \fX_2 \to \G_1$  such that:
\begin{enumerate}
\item  for each $x \in \fX_1$ and $\gamma_1 \in \G_1$, there exists  an open set $x \in U_x \subset \fX_1$ such that 
$ \Phi_2(\alpha(\gamma_1,x))  \circ h | U_x = h \circ \Phi_1(\gamma_1) | U_x$; \\
\item for each $y \in \fX_2$ and $\gamma_2 \in \G_2$, there exists   an open set $y \in V_y \subset \fX_2$ such that 
$ \Phi_1(\beta(\gamma_2,y))  \circ h^{-1} | V_y = h^{-1} \circ \Phi_2(\gamma_2) | V_y$.
\end{enumerate}
The maps $\alpha$ and $\beta$ are not assumed to be cocycles over the respective actions.    
\end{defn}

\begin{remark} \label{rmk-identity}
{\rm 
Suppose that $(\fX_1,\G_1,\Phi_1)$ and $(\fX_2,\G_2,\Phi_2)$ are    actions,   and let  $h \colon \fX_1 \to \fX_2$ be a continuous orbit equivalence.  Form the conjugate action $\Psi_2 \colon \G_2 \times \fX_1 \to \fX_1$ where $\Psi_2 = h^{-1} \circ \Phi_2 \circ h$. It then follows that the identity map is an orbit equivalence between the actions   $(\fX_1,\G_1,\Phi_1)$ and $(\fX_1,\G_2,\Psi_2)$. Thus, we can always reduce to the case where $\fX_1 = \fX_2 = \fX$ and $h$ is the identity map, and if $(\fX,\G_1,\Phi_1)$ is minimal then $(\fX,\G_2,\Phi_2)$ is also minimal.
 }
\end{remark}

 \section{Equicontinuous actions}\label{sec-equicontinuous}

 We  show that equicontinuity is an invariant of continuous orbit equivalence.    The conclusion of Proposition~\ref{prop-equi}, with the stronger assumption that both actions  are  free, was stated in 
Cortez and Medynets         \cite[Corollary~4.4]{CortezMedynets2016}, as a consequence of  Remark~3 in \cite[Section~2]{Medynets2011} that  an isomorphism of full groups is realized spatially for Cantor actions.  
The proof below  follows directly from the definition of a continuous  orbit equivalence.   

   \begin{prop}\label{prop-equi}
 Suppose that    Cantor actions    $(\fX_1,\G_1,\Phi_1)$ and    $(\fX_2,\G_2,\Phi_2)$ are continuously orbit equivalent. If both $\G_1$ and $\G_2$ are finitely generated groups, and   $(\fX_1,\G_1,\Phi_1)$  is   equicontinuous, then so is $(\fX_2,\G_2,\Phi_2)$.
  \end{prop}

  \proof
    By Remark~\ref{rmk-identity}, 
we can assume that  the Cantor spaces are the same, so $\fX = \fX_1 = \fX_2$, and the orbit equivalence $h$ is the identity map on  $\fX$. Let   $\dX$ be   a metric on $\fX$ compatible with the topology. We must show there exists $\e_0 > 0$ so that for any $0 < \e \leq \e_0$ there exists $\delta > 0$ such that for $x,y \in \fX$ with $\dX(x,y) < \delta$,   and   for all $h \in \G_2$ we have  $\dX(\Phi_2(h)(x), \Phi_2(h)(y)) < \e$.
The idea of the proof of this claim  is to show that the action $\Phi_1$ has a  ``shadowing property'', using an idea from the proof of  \cite[Theorem~3.3]{CortezMedynets2016}.

  We   first establish some technical preliminaries.  
    Let $\alpha$ and $\beta$ be the maps in Definition~\ref{def-torb1} for $h$ the identity map. That is, we have continuous maps $\alpha \colon \G_1 \times \fX  \to \G_2$  and $\beta \colon \G_2 \times \fX  \to \G_1$  so that    for $y \in \fX$ and $g \in \G_1$, there exist a clopen set $y \in U_{g,y} \subset \fX$ with
   \begin{equation}\label{eq-conjmap1} 
\Phi_2(\alpha(g,y))(z)      =      \Phi_1(g)(z)   ~ \textrm{for} ~   z \in U_{g,y}   ~ ,
\end{equation}
and for $h \in \G_2$, there exists a clopen set $y \in V_{h,y} \subset \fX$ so that 
 \begin{equation}\label{eq-conjmap2}
\Phi_1(\beta(h,y))(z)      =      \Phi_2(h)(z)   ~ \textrm{for}  ~ z \in V_{h,y} ~ .  
\end{equation}
Let $\Delta(\G_2) \equiv \{h_1, \ldots, h_{\mu}\} \subset \G_2$ be a symmetric set of generators for $\G_2$.  That is, for  $h_i \in \Delta(\G_2)$,   we have  $h_i^{-1} \in \Delta(\G_2)$ for all $1 \leq i \leq \mu$.

For each $1 \leq j \leq \mu$,  we have an open covering of $\fX$ by the sets $\{ V_{h_j, y} \mid y \in \fX\}$. As $\fX$ is compact there exists a Lebesgue number $\e_j >0$ for the covering.
Then $\e' = \min \{\e_1, \ldots , \e_{\mu}\} > 0$ is a Lebesgue number for all of these coverings.

Given  $x \in \fX$   there exists an adapted neighborhood basis at $x$ for the action $\Phi_1$ as in Definition~\ref{def-adaptednbhds}. It follows that we can choose an adapted set  $W \subset \fX$ for the action $\Phi_1$    such that for all $g \in \G_1$, we have $\diam_{\dX} (\Phi_1(g)(W)) < \e'$.
Then the translates $\cW = \{\Phi_1(g)(W) \mid g \in \G_1\}$ form a finite covering of $\fX$ by disjoint clopen sets, and so there is a minimum distance separating them, 
$$\e'' = \min \ \{\dist_{\dX}( \Phi_1(g)(W), \Phi_1(g')(W)) \mid \Phi_1(g)(W) \ne \Phi_1(g')(W)\} > 0 ~ . $$
Then for $0 < \lambda < \e''$ and $y \in \Phi_1(g)(W)$, the ball of radius $\lambda$ about $y$ satisfies $B_{\dX}(y, \lambda) \subset \Phi_1(g)(W)$.

Set $\e_0 = \min \ \{\e', \e''\}$ and choose  $0 < \e < \e_0$.
As the action  $(\fX,\G_1,\Phi_1)$ is equicontinuous, there exists $\delta > 0$ such that  for all $g \in \G_1$ and $x, y \in \fX$ with $ d_X(x,y) < \delta$, then    $\dX(\Phi_1(g)(x), \Phi_1(g)(y)) < \e$. Note that $\delta \leq \e$   as we can take $g$ to be the identity element.  

By the above choices, we have that for $x \in \fX$ and each $1 \leq j \leq \mu$,   there is 
 \begin{itemize}
\item $g_{x} \in \G_1$   such that $x \in  \Phi_1(g_{x})(W)$
\item  $z_{x,j} \in \fX$ so that $B_{\dX}(x, \e) \subset \Phi_1(g_{x})(W) \subset V_{h_j, z_{x,j}}$ 
\end{itemize}
where  the set $V_{h_j, z_{x,j}}$ is defined by \eqref{eq-conjmap2}.

 Now let $x, y \in \fX$ satisfy $\dX(x,y) < \delta$, and let $h \in \G_2$.  We   show that $\dX(\Phi_2(h)(x) , \Phi_2(h)(y)) < \e$.

First,  express $h$ in terms of the generators  $\Delta(\G_2)$, so 
$h = h_{j_m} \cdots h_{j_1}$ for indices $1 \leq j_{\ell} \leq \mu$. We proceed by induction on the factors of $h$.
Set $x_0 = x$, $y_0 = y$, then recursively define for $0 \leq \ell < m$,
 \begin{equation}\label{eq-orbit}
x_{\ell + 1}   =    \Phi_2(h_{j_{\ell + 1}})(x_{\ell}) \quad , \quad y_{\ell+1}   =    \Phi_2(h_{j_{\ell + 1}})(y_{\ell}) \ .
\end{equation}

  Let  $g_{x,0} \in \G_1$ be such that $x_0 \in  \Phi_1(g_{x,0})(W)$, then   we also have $y_0 \in \Phi_1(g_{x,0})(W)$  by the choice of $\delta$.
  Then  there exists $z_0 \in \fX$ such that $B_{\dX}(x_0, \e) \subset  V_{h_{j_1}, z_0}$ and so also $y_0 \in V_{h_{j_1}, z_0}$.
It follows that $\beta(h_{j_1}, x_0) = \beta(h_{j_1}, y_0)  \in \G_1$. 
Set $g_{j_1} = \beta(h_{j_1}, x_0)$ then by \eqref{eq-conjmap2} we have 
$$x_1 = \Phi_2(h_{j_1})(x_0) = \Phi_1(g_{j_1})(x_0) \quad , \quad y_1 = \Phi_2(h_{j_1})(y_0) = \Phi_1(g_{j_1})(y_0) \ .$$
Note that $\dX(x_1 , y_1) < \e$ by the the choice of $\delta$ and the equicontinuity hypothesis for $\Phi_1(g_{j_1})$.

Now let $1 \leq \ell < m$, and assume that  $\{g_{j_1} , g_{j_2}, \ldots , g_{j_{\ell}}\} \subset \G_1$ have been chosen so that for $1 \leq i \leq \ell$ we have  $\dX(x_i , y_i) < \e$ with  
 $$x_{i} = \Phi_2(h_{j_i})(x_{i-1}) = \Phi_1(g_{j_i})(x_{i-1}) \quad , \quad y_i = \Phi_2(h_{j_i})(y_{i-1}) = \Phi_1(g_{j_i})(y_{i-1}) \ .$$
 Then there exists $z_{\ell} \in \fX$ such that 
$B_{\dX}(x_{\ell}, \e) \subset  V_{h_{j_{\ell +1}}, z_{\ell}}$ and so also $y_{\ell} \in V_{h_{j_{\ell +1}}, z_{\ell}}$.

It follows that $\beta(h_{j_{\ell+1}}, x_{\ell}) = \beta(h_{j_{\ell +1}}, y_{\ell})  \in \G_1$. Then set $g_{j_{\ell+1}} = \beta(h_{j_{\ell+1}}, x_{\ell})$ and define $\ds \whg_{\ell+1} \equiv g_{j_{\ell +1}} \cdot g_{j_{\ell}} \cdot \cdots g_{j_{1}}$.
Then   by \eqref{eq-conjmap2} and the previous choices, we have 
\begin{eqnarray*}
x_{\ell + 1} & = &  \Phi_2(h_{j_{\ell + 1}})(x_{\ell}) = \Phi_1(g_{j_{\ell + 1}})(x_{\ell}) = \Phi_1(\whg_{\ell+1})(x_0)  \\
y_{\ell+1} & = &  \Phi_2(h_{j_{\ell + 1}})(y_{\ell}) = \Phi_1(g_{j_{\ell +1}})(y_{\ell}) = \Phi_1(\whg_{\ell+1})(y_0)\ .
\end{eqnarray*}
Then $\dX(x_{\ell +1} , y_{\ell +1}) < \e$ by the the choice of $\delta$  and the equicontinuity hypothesis   for   $\Phi_1(\whg_{\ell+1})$. 
Thus, for   $\ell = m$   we obtain the estimate 
$$\dX(\Phi_2(h)(x) , \Phi_2(h)(y)) = \dX(\Phi_1(\whg_{\mu})(x_0), \Phi_1(\whg_{\mu})(y_0)) < \e$$
as was to be shown. 
  \endproof

 \section{Return equivalence}\label{sec-return}

In this section, we show that the locally quasi-analytic  property of an equicontinuous Cantor action   is preserved by continuous orbit equivalence. The strategy of the proof is to first show that the actions are return equivalent, as defined    in Definition~\ref{def-return}.  Then   Corollary~\ref{cor-LQA} deduces the   locally quasi-analytic property from return equivalence.

In a previous work \cite{HL2019b}, the authors showed   that the  stable property for an equicontinuous  action  is preserved by continuous orbit equivalence. The stable property implies the locally quasi-analytic property, which yields the conclusion of Theorem~\ref{thm-LQA} below for stable actions. However, there exists nilpotent Cantor actions which are not stable \cite{HL2020b}, so we must prove a variant of Theorem~6.10 in \cite{HL2019b}  for our purposes. The proofs of the stable and the locally quasi-analytic versions of this result have significant overlaps, so when possible we   refer to  proofs of the corresponding results  in \cite{HL2019b}.

   \begin{thm}\label{thm-LQA}
Let    $(\fX_1,\G_1,\Phi_1)$ and    $(\fX_2,\G_2,\Phi_2)$ be   Cantor actions, with  both $\G_1$ and $\G_2$   finitely generated groups. Suppose that the actions are continuously orbit equivalent, and that 
   $(\fX_2,\G_2,\Phi_2)$  is   equicontinuous and  locally quasi-analytic. Then   $(\fX_1,\G_1,\Phi_1)$ is return equivalent to $(\fX_2,\G_2,\Phi_2)$.
  \end{thm}
  
 \proof
 
    By Remark~\ref{rmk-identity}, 
we can assume that  the Cantor spaces are the same, so $\fX = \fX_1 = \fX_2$, and the orbit equivalence $h$ is the identity map on  $\fX$. Let   $\dX$ be   a metric on $\fX$ compatible with the topology. 
Let   $\alpha \colon \G_1 \times \fX  \to \G_2$  and $\beta \colon \G_2 \times \fX  \to \G_1$    satisfy the relations \eqref{eq-conjmap1}  and \eqref{eq-conjmap2}. Note that  $(\fX_1,\G_1,\Phi_1)$ is an equicontinuous action by Proposition~\ref{prop-equi}. 
 
The proof  that the actions are return equivalent    follows from a sequence of results.  We first show in Lemma~\ref{lem-cocycle} that there exists an adapted set $U \subset \fX$ such that the map   $\alpha$ restricted to the action of $\G_{1,U}$ on $U$ satisfies the cocycle identity. 
 We then show in Proposition~\ref{prop-coboundary} that there exists an adapted set $W \subset U$ such that the cocycle $\alpha$ is a coboundary when restricted to the action of $\G_{1,W}$ on $W$, and thus conclude   that it induces  a group homomorphism on $\G_{1,W}$. Both of these results have their exact counterparts in the proof of \cite[Theorem~6.10]{HL2019b}, so we only include sufficient details of their proofs to establish the notation required for the proof of the new results, 
  Lemmas~\ref{lem-inducedhomo}, \ref{lem-adapted} and \ref{lem-returned}, which  show  that this  homomorphism induces a return equivalence of the actions.

 Choose $x_0 \in \fX$, then as  $\Phi_2$ is locally quasi-analytic, there exists $V \subset \fX$ adapted to the action $\Phi_2$   such that  $x_0 \in V$ and  the action of  $H_{2,V} = \Phi_2(\G_{2,V}) \subset \Homeo(V)$ on $V$ is topologically free. Thus, there exists a dense   subset $\cZ_V \subset V$ such that the action of $H_{2,V} $ is free when restricted to $\cZ_V$.

Choose    an adapted set $U \subset \fX$   for the action $\Phi_1$ with $x_0 \in U \subset V$.

Let $\G_{1,U} \subset \G_1$ be the isotropy group of the action $\Phi_1$ on $U$, as defined as in  \eqref{eq-adapted}. 

  Let $\alpha_U \colon  \G_{1,U} \times U \to \G_2$ denote the restriction of   the map $\alpha \colon \G_1 \times \fX  \to \G_2$.  
  Then for  each $g \in \G_{1,U}$ and $y \in U$ we have $\Phi_1(g)(y) \in U$. 
 Set $h = \alpha(g,x_0) \in \G_2$ so that  by \eqref{eq-conjmap1} we have  $\Phi_2(h)(y) = \Phi_1(g)(y)$. 
 
 Then $U \subset V$ implies that $\Phi_2(h)(V) \cap V \ne \emptyset$ hence $h \in \G_{2,V}$ and 
  $\Phi_{2,V}(h) \in H_{2,V}   \subset \Homeo(V)$. 
  
  That is,  
  the restriction of $\alpha$ to $\G_{1,U} \times U$ induces a map  $\whalpha_U = \Phi_{2,V} \circ \alpha \colon  \G_{1,U} \times U \to H_{2,V} $.
The action of $H_{2,V} $ on $V$ is topologically free,  so we have:

 \begin{lemma}\label{lem-cocycle}
 $\whalpha_U \colon  \G_{1,U} \times U \to H_{2,V} $ satisfies the cocycle identity 
 \begin{equation}\label{eq-cocycleidentity}
\whalpha_U(g \cdot g', y) = \whalpha_U(g, \Phi_1(g')(y)) \cdot \whalpha_U(g',y) \quad \text{ for } g,g' \in \G_{1,U} \text { and } y \in U \ .
\end{equation}
 \end{lemma}
   \proof
  This follows as in the proof of   \cite[Lemma~2.8]{Li2018}, or that of      \cite[Proposition~6.12]{HL2019b}.   
   \endproof

The next result asserts that a properly chosen restriction of the cocycle $\whalpha_U$ in Lemma~\ref{lem-cocycle} is a coboundary. The proof uses in an essential way that the group $\G_1$ is finitely-generated, as this allows factoring the cocycle  $\whalpha_U$   into actions supported on regions of continuity for the continuous orbit equivalence functions in 
  \eqref{eq-conjmap1} and \eqref{eq-conjmap2}.  The idea of the proof  is modeled on that   of \cite[Theorem~3.3]{CortezMedynets2016}, with the variation   that we   only assume   the   action of $H_{2,V} $ is topologically free, and do not assume   the   action $(\fX_1,\G_1,\Phi_1)$   is locally quasi-analytic.

\begin{prop}\label{prop-coboundary}
Let    $(\fX_1,\G_1,\Phi_1)$ and    $(\fX_2,\G_2,\Phi_2)$ be equicontinuous Cantor actions, with  both $\G_1$ and $\G_2$   finitely generated groups. Suppose that the actions are continuously orbit equivalent, and that 
   $(\fX_2,\G_2,\Phi_2)$  is    locally quasi-analytic. Then there exists an adapted set $W \subset \fX$ for the action $\Phi_1$  so that the restricted cocycle 
    $\whalpha_W \colon  \G_{1,W} \times W \to H_{2,V} $ is   induced by a  group homomorphism 
$\whtheta_W \colon \G_{1,W} \to H_{2,V} $. That is, for $g \in \G_{1,W}$  and $x \in W$,  we have 
$\whalpha_W(g,x) = \whtheta_W(g)$.
 \end{prop}  
 \proof
  This follows exactly as in the proof of \cite[Proposition~6.12]{HL2019b}. 
  \endproof

 The next steps in the proof of Theorem~\ref{thm-LQA}   deviate from that of \cite[Theorem~6.10]{HL2019b}, as we must show that there exists an isomorphism of holonomy actions as in Definition~\ref{def-return}. This follows from the next three results.   
 Set     
  $H_{1,W} = \Phi_{1_W}(\G_{1,W}) \subset \Homeo(W)$, and $H_{2,V} = \Phi_2(\G_{2,V}) \subset \Homeo(V)$. 

 \begin{lemma}\label{lem-inducedhomo}
   $  \whtheta_W$  induces a monomorphism  $\theta_W \colon H_{1,W} \to H_{2,V}$. 
 \end{lemma}
 \proof
 We first show that 
 \begin{equation}\label{eq-quotientmap}
\ker \{\Phi_{1,W} \colon \G_{1,W} \to H_{1,W} \} \subset \ker \{\whtheta_W \colon \G_{1,W} \to H_{2,V} \} \ .
\end{equation}
 Suppose  that   $g \in \G_{1,W}$ satisfies $\Phi_{1,W}(g) = \Id \in H_{1,W} \subset \Homeo(W)$. 
 Recall that for  $y \in \cZ_V \cap W$  the action of $H_{2,V}$ is free on the orbit of $y$, and  
 we  have that $\Phi_{1,W}(g)(y) = y$. For $h = \whtheta_W(g)$, 
  by  the   identity  \eqref{eq-conjmap1}  we have
  \begin{equation}
h \cdot y = \whtheta_W(g)(y)  = \whalpha_W(g,y) \cdot y = \Phi_{1,W}(g)(y) = y \ .
\end{equation}
Thus, $h \in H_{2,W}$ must be the identity map since $y \in \cZ_V$, so $\whtheta_W(g) = \Id$.  Thus \eqref{eq-quotientmap} is satisfied, so  $\whtheta_W$  induces   a well-defined homomorphism  $\theta_W \colon H_{1,W} \to H_{2,V}$.

 Suppose that $\theta_W(\whg) = \Id$ for $\whg \in H_{1,W}$, then for $g \in \G_{1,W}$ with $\whg = \Phi_{1,W}(g)$, using the identity  \eqref{eq-quotientmap} again, we have that $h \in H_{2,W}$ is the identity map,  so $\theta_W$ is injective.
 \endproof

  \eject
  
  \begin{lemma}\label{lem-adapted}
 The adapted set $W \subset \fX$  for the action $\Phi_1$  is also adapted for the action of $\Phi_2$.
 \end{lemma}
\proof
Let $h \in \G_2$ be such that $\Phi_2(h)(W) \cap W \ne \emptyset$. As $W \subset V$ is a clopen set and $\cZ_V \subset V$ is dense,  the set $\cZ_V \cap W$ is dense in $W$. As $W  \cap \Phi_2(h^{-1})(W)$ is non-empty and open in $W$, there exists     $y\in \cZ_V \cap W  \cap \Phi_2(h^{-1})(W)$  for which we have $\Phi_2(h)(y) \in W$. 

Set $g = \beta(h,y) \in \G_1$ then by  \eqref{eq-conjmap2}    we have 
$\Phi_1(g)(y)   = \Phi_2(h)(y)    \in W$.
As $W$ is adapted to the action $\Phi_1$ we have $\Phi_1(g)(W) = W$ so $g \in \G_{1,W}$.

Now set $h' = \whalpha_W(g,x) = \whtheta_W(g) \in \G_{2,V}$, where $\whtheta_W$ is the map defined in
Proposition~\ref{prop-coboundary}. Then $\Phi_{2,V}(h')(y) = \Phi_{2,V}(h)(y)$ and so $\Phi_{2,V}(h') = \Phi_{2,V}(h) \in H_{2,V}$ as $y \in \cZ_V$. That is, $h = \whtheta_W(g)$ and so for all $z \in W$ we have $\Phi_{2,V}(h)(z) = \Phi_{1,W}(g)(z) \in W$.

Thus,   $W$ is adapted to the action $\Phi_2$ as was to be shown.
\endproof
  
  
 \begin{lemma}\label{lem-returned}
 For the adapted set $W$, the map $\whtheta_W$ induces an isomorphism $\theta_W \colon H_{1,W} \to H_{2,W}$.
 \end{lemma}
\proof
For $g \in \G_{1,W}$ set $\whg = \Phi_{1,W}(g) \in H_{1,W}$. The proofs of Lemmas~\ref{lem-inducedhomo} and \ref{lem-adapted} show that $\theta_W(\whg) \in H_{2,W}$. Given $h \in \G_{2,W}$ and $y \in \cZ_V \cap W$ set $g = \beta(h,y) \in \G_1$. Then $\whtheta_W(g)$ and $\Phi_2(h)$ agree on an open neighborhood of $y$ in $W$, hence agree on all of $W$ as $W \subset V$. It follows that $\theta_W$ is an isomorphism onto.
\endproof
  
We have shown that  $W \subset \fX$ is adapted to both actions $\Phi_1$ and $\Phi_2$,  and    $\theta_W \colon H_{1,W} \to H_{2,W} \subset H_{2,V}$ is an isomorphism.  This completes the proof of Theorem~\ref{thm-LQA}.
\endproof

   \begin{cor}\label{cor-LQA}
Let    $(\fX,\G_1,\Phi_1)$ and    $(\fX,\G_2,\Phi_2)$ be   Cantor actions, with  both $\G_1$ and $\G_2$   finitely generated groups. Suppose that the identity map on $\fX$ is a   continuous orbit equivalence, and that 
   $(\fX,\G_2,\Phi_2)$  is   equicontinuous and  locally quasi-analytic. Then   $(\fX,\G_1,\Phi_1)$ is     locally quasi-analytic.
  \end{cor}
\proof
It follows from Theorem~\ref{thm-LQA}  that the two actions are return equivalent, for an adapted set $W \subset \fX$. As $(\fX,\G_2,\Phi_2)$ is locally quasi-analytic, we can chose $W$ sufficiently small so that the induced action of $H_{2,W}$ on $W$   is topologically free. Then  the isomorphic action of $H_{1,W}$ on $W$ is also topologically free, and thus the action of $\G_{1,W}$ on $W$ is quasi-analytic. As $W$ is adapted for the action $\Phi_1$, it follows that $(\fX,\G_1,\Phi_1)$ is locally quasi-analytic.
\endproof

  \section{Nilpotent actions} \label{sec-nilpotent}
 
 In this section, we give the proof of Theorem~\ref{thm-nilconjugate}, and also give   examples to illustrate its conclusions.
  
Let $(\fX_1,\G_1,\Phi_1)$ be a nilpotent Cantor action. The group $\G_1$ satisfies the Noetherian property \cite{Baer1956} 
for increasing chains of subgroups, so  the action is locally quasi-analytic by the following result:
\begin{thm}\cite[Theorem~1.6]{HL2018b}
Let    $\G$ be a Noetherian group. Then   a  minimal equicontinuous Cantor action $(\fX,\G,\Phi)$   is locally quasi-analytic.
\end{thm}

Let $(\fX_2,\G_2,\Phi_2)$ be a Cantor action, and assume that $\G_2$ is finitely-generated.

Assume that the actions are continuously orbit equivalent.      By Remark~\ref{rmk-identity}, 
we can assume that  the Cantor spaces are the same, so $\fX = \fX_1 = \fX_2$, and the orbit equivalence $h$ is the identity map.

Then by Proposition~\ref{prop-equi}, the action $(\fX_2,\G_2,\Phi_2)$ is equicontinuous, and by  
  Theorem~\ref{thm-LQA}, the actions  $(\fX_1,\G_1,\Phi_1)$ and $(\fX_2,\G_2,\Phi_2)$ are return equivalent. Then by Corollary~\ref{cor-LQA} the action  $(\fX_2,\G_2,\Phi_2)$   is locally quasi-analytic. 
  Let  $W \subset \fX$ be the clopen set adapted to both actions   $\Phi_1$ and $\Phi_2$, chosen as in  the proof of Corollary~\ref{cor-LQA} so that  both  actions restricted to $W$ are quasi-analytic.

Let $\G_{1,W} \subset \G_1$ be the isotropy subgroup of $W$ for the action $\Phi_1$,   with holonomy group $H_{1,W} = \Phi_{1,W}(\G_{1,W}) \subset \Homeo(U)$. Let $\G_{2,W} \subset \G_2$ be the isotropy subgroup of $W$ for the action $\Phi_2$, with holonomy group $H_{2,W} = \Phi_{2,W}(\G_{2,W}) \subset \Homeo(W)$.

Let  $\theta_W \colon H_{1,W} \to H_{2,W}$ be  the   isomorphism defined in Lemma~\ref{lem-returned} which implements the orbit equivalence between the two actions.
As $\G_{1,W} \subset \G_1$ has finite index, there exists a nilpotent subgroup $\Lambda_1 \subset \G_{1,W}$ of finite index, with $\Lambda_1$ finitely generated.
Then the image $\Lambda_2 = \whtheta_W(\Lambda_1) \subset H_{2,W}$ is a finitely-generated nilpotent subgroup  of finite index.

Suppose that the action $\Phi_1$ is topologically free, then the restriction $\Phi_{1,W} \colon \G_{1,W} \to H_{1,W}$ is an isomorphism. Likewise, if   $\Phi_2$ is topologically free, then the restriction $\Phi_{2,W} \colon \G_{2,W} \to H_{2,W}$ is an isomorphism.  As $\Lambda_1 \subset \G_{1,W}$  has finite index, and likewise for $\Lambda_2 \subset \G_{2,W}$, this shows that the groups $\G_{1,W}$ and $\G_{2,W}$   contain isomorphic nilpotent subgroups of finite index, and thus also  $\G_1$ and $\G_2$. In particular, the groups $\G_1$ and $\G_2$   are commensurable.

Now, assume  that the action $(\fX_2,\G_2,\Phi_2)$  is effective, that is,  the  action map $\Phi_2 \colon \G_2 \to \Homeo(\fX)$ is injective.
Let $K_{2,W} \subset \G_{2,W}$ be the kernel of  the restricted map 
$\Phi_{2,W} \colon \G_{2,W} \to H_{2,W}$, which need not be trivial (see Examples~\ref{ex-comm} and \ref{ex-notcomm} below).  
  
 Let $X_{2,W} = \G_2/\G_{2,W}$ be the finite set of cosets of $\G_{2,W}$, with a transitive left $\G_2$ action.  The action $\Phi_2$ induces a map $\Pi_{2,W} \colon \G_2 \to \Perm(X_{2,W})$ to the group of permutations of $X_{2,W}$, and $\G_{2,W}$ is the isotropy subgroup  of the identity coset $e_W \in  X_{2,W}$. 
  Let $C_{2,W} = \ker (\Pi_{2,W}) \subset \G_2$ be the kernel of this representation, so $C_{2,W}$ is a normal subgroup of $\G_2$ with finite index.

Choose representatives $\{h_i \in \G_2 \mid   1 \leq i \leq \nu\}$ of the cosets of $\G_2/\G_{2,W}$ and set $W_i = h_{i} \cdot W$.
Then 
$$\fX = W_1 \ \cup  W_2 \ \cup \ \cdots \ \cup \  W_{\nu} \ .$$ 
For  $h \in C_{2,W}$, the action of $\Phi_2(h)$  on $\fX$ leaves each clopen set $W_i$ invariant, so for $y =  h_i \cdot z\in W_i$ with $ z \in W$, we have:
\begin{equation}\label{eq-conjugaterel}
h \cdot y = h \cdot h_i \cdot z  = h_i \cdot (h_i^{-1} \ h  \ h_i) \cdot z   
\end{equation}
where $h_i^{-1} h \ h_i \in C_{2,W}$,  as $C_{2,W}$ is normal in $\G_2$.
For $1 \leq i \leq \nu$,  define the conjugate action on $\fX$,   
 $$\Phi_2^i(h)(z) = \Phi_2(h_i)^{-1} \Phi_2( h) \Phi_2(h_i)(z) \ , \quad z \in \fX \ .$$
 Then  for  $h \in C_{2,W}$, by \eqref{eq-conjugaterel} the restriction $\Phi_2(h)  \colon  W_i \to W_i$ is the identity if and only if 
 $$h \in \ker \{\Phi_2^i \colon C_{2,W} \to H_{2,W} \subset \Homeo(W)\} \ .$$

For $h \in C_{2,W}$ which is not the identity, we have by assumption that $\Phi_2(h)$ is not the identity map on $\fX$, hence there exists some $1 \leq i \leq \nu$ such that $\Phi_2(h)   \colon  W_i \to  W_i$ is not the identity, and so $h \not\in \ker \{\Phi_2^i \colon C_{2,W} \to H_{2,W} \}$.

Define a representation $\whrho_2$ of $C_{2,W}$ into a product of $\nu$ copies of $H_{2,W}$ by setting for $h \in C_{2,W}$, 
\begin{equation}\label{eq-prodrep}
\whrho_2 \colon C_{2,W} \to H_{2,W} \times \cdots \times H_{2,W} \quad , \quad \whrho_2(h) = \Phi_2^1(h)  \times \cdots \times \Phi_2^{\nu}(h) 
\end{equation}
The kernel of  $\whrho_2$ is trivial by the above arguments  and the assumption that the action $\Phi_2$ is effective.

Recall that $\Lambda_2 = \theta_W(\Lambda_1)\subset H_{2,W}$ is a finitely-generated nilpotent subgroup  of finite index.

 For $1 \leq i \leq \nu$, let $\Lambda_2^i = (\Phi_2^i)^{-1}(\Lambda_2)$. Then $\Lambda_2^i$ is a subgroup of finite index in $C_{2,W}$, and so $\Lambda_2' = \Lambda_2^1 \cap \cdots \cap \Lambda_2^{\nu}$ has finite index in $C_{2,W}$ and thus also in $\G_2$.
 Observe that for each $1 \leq i \leq \nu$, we have $ \Phi_2^i(\Lambda_2') \subset \Lambda_2$. Moreover, the homomorphism \eqref{eq-prodrep} restricts to an embedding 
 \begin{equation}\label{eq-prodrep2}
\whrho_2 \colon \Lambda_2' \to \Lambda_2 \times \cdots \times \Lambda_2 \ .
\end{equation}

Thus, $\whrho_2$ is an injection of $\Lambda_2'$ into a product of nilpotent groups, which is again nilpotent, and so $\Lambda_2'$ is a nilpotent group.
Hence, $\G_2$ is virtually nilpotent, as was to be shown.
 \endproof

 \begin{ex}\label{ex-comm}
{\rm 
 We give an elementary example to show that the conclusion  that the groups $\G_1$ and $\G_2$ are commensurable in Theorem~\ref{thm-nilconjugate} is best possible. 
 
 Let $\G_0$ be a finitely-generated, torsion free, infinite nilpotent group, and $\fX_0 = \whGamma_0$ be the profinite completion of $\G_0$.  Let $\Phi_0 \colon \G_0 \times \fX_0 \to \fX_0$ be the action by left multiplication. Then the action $\Phi_0$  on $\fX_0$ is free.

 Let $Q_1 = \mZ/2\mZ \times \mZ/2\mZ$ be the product of cyclic groups of order two, and let $Q_2 = \mZ/4\mZ$ be the cyclic group of order 4. Let $Y= \{0,1,2,3\}$ be the set with 4 elements.  Choose  identifications $\tau_1 \colon Q_1 \to Y$ and $\tau_2 \colon Q_2 \to Y$, which define actions of $Q_1$ and $Q_2$ on $Y$. 
 
   Let $\fX = \fX_0 \times Y$ be the product Cantor space. Define $\G_1 = \G_0 \times Q_1$ with the product action $\Phi_1$ on  $\fX$. Similarly, let  $\G_2 = \G_0 \times Q_2$ with the product action $\Phi_2$ on $\fX$. Both actions are minimal, equicontinuous and free.   Moreover, both actions have the same orbits on $\fX$ and the orbit map satisfies the conditions in Definition~\ref{def-torb1}. However, $\G_1$ and $\G_2$ are not isomorphic as their characteristic torsion subgroups are not isomorphic. On the other hand, both contain the subgroup $\G_0$ with finite index, so are commensurable.
 }
 \end{ex}

 \begin{ex}\label{ex-notcomm}
\rm{
  We give an example to show that the hypothesis that the actions $\Phi_1$ and $\Phi_2$ are topologically free  in Theorem~\ref{thm-nilconjugate} is necessary to conclude that $\G_1$ and $\G_2$ are commensurable.

  Let $\G_0$ be a finitely-generated, torsion free, infinite nilpotent group, and $\fX_0 = \whGamma_0$ be the profinite completion of $\G_0$.   
Let $\Phi_0 \colon \G_0 \times \fX_0 \to \fX_0$ be the free action by left multiplication.

Choose a non-trivial  finite group $Q_0$.     
 List the elements $Q_0 = \{q_1, \ldots, q_k\}$ with $q_1$ the identity element, then write $\fX = Q_0 \times \fX_0$ as the union of clopen subsets $W_{i} =  \{q_i\} \times \fX_0$, so 
 $\fX = W_1 \cup \cdots \cup W_k$. 
 Let $Q_0$ act on $\fX$ as the identity on the factor $\fX_0$, and 
 by left multiplication on the factor $Q_0$, so that it transitively permutes the partition $\{W_1, \ldots , W_k\}$ of $\fX$.  

We now define two minimal equicontinuous actions $\Phi_1$ and $\Phi_2$   on $\fX$, where $\Phi_1$ is a free action,    $\Phi_2$ is locally quasi-analytic, and the actions are continuously orbit equivalent, but the  groups $\G_1$ and $\G_2$ are not commensurable.

First, define  $\G_1 = Q_0 \times \G_0$ and the action $\Phi_1$ is defined as follows. The action of $Q_0$ on $\fX$ is that above,  
while the $\Phi_1$-action of $\G_0$ acts as the identity on the set $Q_0$ and as   the     $\Phi_0$-action of $\G_0$ on $\fX_0$.
Note that this action is free.

  Next, define    $\ds \G_2 =  Q_0  \ltimes \ \G_0 $ as the wreath product, namely, let $\G_0^{|Q_0|} = \{f: Q_0 \to \G_0\}$ be the set of functions, and note that there is a shift action
  $$Q_0 \times \G_0^{|Q_0|} \to \G_0^{|Q_0|} : (q,f) \mapsto f^q (*) = f(q^{-1}   (*)).$$
Then the \emph{wreath product} $Q_0 \ltimes \G_0^{|Q_0|} $
  is a group with group product $$(q_1,f_1)\circ (q_2,f_2)=(q_1 q_2, f_1 f_2^{q_1}).$$ 
  The wreath product $\G_2$ acts on $Q_0 \times \fX_0$ by
    \begin{eqnarray}\label{wp-action} (q,f)\cdot (s,x) = (qs, f(qs)(x)). \end{eqnarray}
That is, the action \eqref{wp-action} permutes the copies of $\fX_0$ in the product $Q_0\times \fX_0$, while acting on each copy of $\fX_0$ independently via an element defined by the function $f$.

 We show that the identity map on $\fX$ is a continuous orbit equivalence between the actions   $(\fX, \G_1, \Phi_1)$ and $(\fX, \G_2, \Phi_2)$. For that, we define cocycles $\alpha: \G_1 \times (Q_0 \times \fX_0) \to \G_2$ and $\beta: \G_2 \times (Q_0 \times \fX_0) \to \G_1$ as in Definition \ref{def-torb1}. 
 
Let $f_g: Q_0 \to \G_0^{|Q_0|}: q \to g$ be the constant function, and define the following function, which is independent of the second component, and so it satisfies (1) in Definition \ref{def-torb1},
  $$\alpha((q,g), (s,x)) \mapsto (q,f_g).$$
 The function $\alpha$ implements a ``diagonal''  embedding of $\G_1 \subset \G_2$. A straightforward computation using \eqref{wp-action} shows that  the orbits of 
 $(\fX, \G_1, \Phi_1)$ are contained in the orbits of $(\fX, \G_2, \Phi_2)$. 
 
 Conversely, the following function is independent of $x$ and so it is constant on the clopen sets $W_i$, for $i = 1,\ldots, k$,
   $$\beta((q,f)(s,x)) = (q,f(qs)).$$
 Thus $\beta$ satisfies (2) in Definition \ref{def-torb1}, and clearly maps orbits of $(\fX, \G_2, \Phi_2)$ to orbits of $(\fX, \G_1, \Phi_1)$. Thus $\alpha$ and $\beta$ implement a continuous orbit equivalence between the two actions.

Set $U = W_1$. Then  $\G_{1,U} = \{e\} \times \G_0$. On the other hand, 
$\ds \G_{2,U} = \{e\} \times \G_0^{|Q_0|}  \cong  \{e\} \times \oplus_{i=1}^k \ \G_0^i $, where $\G_0^i = \G_0$.
Thus, the  groups $\G_1$ and $\G_2$ are not commensurable.

 }
 \end{ex}

\begin{remark}
{\rm 
The idea of Example~\ref{ex-notcomm} is that while an orbit equivalence between   actions fixes their orbits,   for locally quasi-analytic actions  it does not determine   the actions of the isotropy groups of clopen sets. This is seen in the above example where  the   isotropy groups $\G_{1,U}$ and $\G_{2,U}$ are related, but not isomorphic.
This construction admits various generalizations. It should also be compared with the proof of Theorem~\ref{thm-vc} below.}
\end{remark}

 
\section{Virtual nilpotency class}\label{sec-class}

In this section we introduce  a property for finitely-generated virtually nilpotent groups which is used to define an invariant for nilpotent Cantor actions.  Let $\Lambda$ be a finitely-generated torsion-free nilpotent group. 
 The \emph{nilpotency class} $\textrm{c}(\Lambda)$ is the least integer $k$ such that for the lower central series 
$\Lambda_0 = \Lambda$, 
$\Lambda_{i+1} = [\Lambda, \Lambda_{i}]$, we have $\Lambda_{k+1} = \{e\}$. 
Note that  $\textrm{c}(\Lambda') = \textrm{c}(\Lambda)$    for any subgroup of finite index $\Lambda' \subset \Lambda$.

 \begin{lemma}\label{lem-tfsubgroup}
Let    $\Lambda$ be a finitely-generated nilpotent group.
Then there exists a finitely-generated torsion-free subgroup $\Lambda' \subset \Lambda$ of finite index.
\end{lemma}
\proof
A finitely generated nilpotent group is  residually finite, hence there exists a descending chain of finite index normal subgroups $\{\Lambda_{\ell} \mid \ell \geq 0\}$ where $\Lambda_0 = \Lambda$ and $\cap_{\ell > 0} \ \Lambda_{\ell} = \{0\}$. 
 Let $\Lambda_t \subset \Lambda$ be the maximal subgroup of torsion elements, then $\Lambda_t$ is finitely generated, hence is a finite group.  Moreover, $\Lambda_t$ is normal in $\Lambda$, and  $\Lambda_t$ contains every element of finite order in $\Lambda$.  It follows that there exists $\ell_0 > 0$ such that $\Lambda_t \cap \Lambda_{\ell_0} = \{0\}$. Then set $\Lambda' = \Lambda_{\ell_0}$.
\endproof

Now let $\G$ be a virtually nilpotent group, so there exists a finitely generated nilpotent subgroup $\Lambda \subset \G$ of finite index. 
Then by     Lemma~\ref{lem-tfsubgroup} there exists a torsion-free subgroup $\Lambda' \subset \Lambda$ of finite index, so we can assume without loss of generality that $\Lambda$ is torsion free. Moreover,  the value $\textrm{c}(\Lambda)$ is independent of the choice of such $\Lambda'$ by the previous remarks.

\begin{defn}\label{def-vnc}
Let   $\G$ be a virtually nilpotent group.   The \emph{virtual nilpotency} class $\vc(\G)=\textrm{c}(\Lambda)$ where $\Lambda \subset \G$ is a torsion-free nilpotent subgroup of finite index.
\end{defn}

Observe that   $\vc(\G) = 0$ implies that $\G$ is a finite group, and   $\vc(\G) = 1$ implies that  $\G$ contains an abelian subgroup of finite index. The discrete Heisenberg group $\cH$ has $\vc(\cH) = 2$.  In addition, note that there are many torsion-free nilpotent groups $\Lambda$ with $c(\Lambda) =2$ that are not congruent to the Heisenberg group (see \cite{LL2002} for example.)

 If $\G_1$ and $\G_2$ are virtually nilpotent groups which are commensurable, that is, have subgroups of finite index which are isomorphic, then $\vc(\G_1) = \vc(\G_2)$.

The following   property of the virtual nilpotency class  will be used to prove Theorem~\ref{thm-vc}.

\begin{lemma}\label{lem-multiplicative}
Let $\Lambda$ be a finitely-generated nilpotent group, and $k \geq 1$ an arbitrary integer. 
Then for the product group $\Lambda^k =  \prod_{i=1}^k \ \Lambda_i$ where each $\Lambda_i = \Lambda$, we then have   $\vc(\Lambda) = \vc \left( \Lambda^k \right)$.
\end{lemma}
 
\proof
Let $\Lambda_0 \subset \Lambda$ be a subgroup of finite index with $c(\Lambda_0) = \vc(\Lambda)$.
Then $$\vc \left( \Lambda^k \right) \leq c\left(  \Lambda_0^k \right) = c(\Lambda_0) = \vc(\Lambda) \ .$$
Conversely, let $D \subset   \Lambda^k$ have finite index with $c(D) = \vc \left(   \Lambda^k \right)$. 

Then $\Lambda_0 = D  \cap \left( \Lambda \times \prod_{i=2}^k \ \{e\} \right)$ has finite index in $\Lambda$ and satisfies 
$$\vc(\Lambda) \leq c(\Lambda_0) \leq c(D) \leq \vc \left(   \Lambda^k \right)$$
which shows the claim.
\endproof

We now give the proof of Theorem~\ref{thm-vc}.
Let   $(\fX_1,\G_1,\Phi_1)$ and    $(\fX_2,\G_2,\Phi_2)$ be effective Cantor actions, with  both $\G_1$ and $\G_2$   finitely generated groups, and  assume that    the actions are continuously orbit equivalent.
If $(\fX_1,\G_1,\Phi_1)$ is a nilpotent Cantor action, then by Theorem~\ref{thm-nilconjugate},  the action  $(\fX_2,\G_2,\Phi_2)$ is   return equivalent to $(\fX_1,\G_1,\Phi_1)$.  We do not assume that the action are topologically free, therefore, our statement does not follow directly from the second sentence in Theorem \ref{thm-nilconjugate}. The statement which we prove here in Theorem \ref{thm-vc} is weaker. It proves that the groups $\G_1$ and $\G_2$ have the same virtual nilpotency class, but this need not imply that they contain isomorphic subgroups.

 Without loss of generality we may assume that $\fX_1 = \fX_2 = \fX$ and that the identity map is an orbit equivalence. Then by the proof of Theorem~\ref{thm-LQA}, there exists an adapted set $W \subset \fX$ for both actions, and an isomorphism $\theta_W \colon H_{1,W} \to H_{2,W}$.

 We next proceed as  in the proof of Theorem~\ref{thm-nilconjugate}. 
  Let $k \geq 1$ be the index of $\G_{1,W}$ in $\G_1$. As the action $\Phi_1$ is effective, we have an injective map
  $\whrho_1 \colon C_{1,W} \to \prod_{i=1}^k \ H_{1,W}$ as in   \eqref{eq-prodrep}. Similarly, as the action $\Phi_2$ is effective, we also have an injective map
  $\whrho_2 \colon C_{2,W} \to \prod_{i=1}^k \ H_{2,W}$, with  the same index $k$.  Indeed, $W$ is adapted to both actions, which implies that the index of $\G_{1,W}$ in $\G_1$ equals the index of $\G_{2,W} $ in $\G_2$.

   Let $\Lambda_1 \subset  C_{1,W} \subset \G_{1,W}$ be a nilpotent subgroup of finite index, and without loss of generality we may assume that $c(\Lambda_1) = \vc(\G_{1,W})$. Then the image $\Lambda_1' = \Phi_{1,W}(\Lambda) \subset H_{1,W}$ satisfies 
   $\vc(\Lambda_1') \leq c(\Lambda_1)$.  
      
   On the other hand, as $\whrho_1$ is injective on $C_{1,W}$,  we have  
   $$c(\Lambda_1) = c(\whrho_1(\Lambda_1)) \leq \vc\left(  \prod_{i=1}^k \ H_{1,W}\right) = \vc\left(  \prod_{i=1}^k \ \Lambda_1' \right) = \vc(\Lambda_1') \ .$$
   Thus, $\vc(\G_{1,W}) = c(\Lambda_1) = \vc(\Lambda_1')$.
   
   By an analogous argument, we have $\vc(\G_{2,W}) = c(\Lambda_2) = \vc(\Lambda_2')$. Since $ \theta_W \colon  \Lambda_1' \to \Lambda_2'$ is an isomorphism,  $\vc(\G_1) = \vc(\G_{1,W}) =\vc(\G_{2,W}) = \vc(\G_2)$. This shows the claim of  Theorem~\ref{thm-vc}.


\end{document}